\documentclass{amsart}
\usepackage{amssymb}
\usepackage{upref}

\usepackage[all,cmtip]{xy}

\usepackage[lite,initials,shortalphabetic]{amsrefs}

\makeatletter
  \renewcommand{\p@enumi}{\thesubsection}
\makeatother

\newenvironment{resumeenumerate}[1]
{\begin{enumerate}
 \setcounter{enumi}{#1}
 \addtocounter{enumi}{-1}
}
{\end{enumerate}
}

\newenvironment{lettered}
{\begin{list}{\thelettercounter)}
 {\usecounter{lettercounter}\def\makelabel##1{\hss\llap{##1}}}
}
{\end{list}
}
\newcounter{lettercounter}
\renewcommand{\thelettercounter}{\alph{lettercounter}}

\makeatletter
\newcommand{\emsection}[1]{%
  \par
  \addpenalty\@secpenalty
  \vskip 6 pt plus 9 pt
  \emph{#1.}\nobreak\enspace\ignorespaces
}
\makeatother

\newcommand{\intro}{%
  \goodbreak
  \vskip 6 pt plus 9 pt
}

\numberwithin{equation}{subsection}

\newcommand{\cat}[1]{\boldsymbol{#1}}

\newcommand{\RelCat}{\mathbf{RelCat}}

\newcommand{\Cat}{\mathbf{Cat}}

\newcommand{\SCat}{\mathbf{SCat}}

\newcommand{\RelSCat}{\mathbf{Rel\,SCat}}

\newcommand{\RelPSCat}{\mathbf{Rel}\,\mathrm{P}\,\textbf{SCat}}

\newcommand{\simp}{\mathrm{s}}

\DeclareMathOperator{\Ho}{Ho}

\newcommand{\LH}{L^{H}}

\newcommand{\intersect}{\cap}

\newcommand{\op}{^{\mathrm{op}}}

\newcommand{\spacedcdots}{{\cdot\;\cdot\;\cdot}}

\newcommand{\adj}[4]{#1\negmedspace: #2\adjarrows #3:\negmedspace #4}
\newcommand{\adjarrows}{\mathchoice{\longleftrightarrow}
  {\leftrightarrow}
  {\leftrightarrow}
  {\leftrightarrow}}

\begin{document}

\title[Partial model categories and their simplicial nerves]
  {Partial model categories\\
  and their simplicial nerves}

\author{C. Barwick}
\address{Department of Mathematics, Massachusetts Institute of
  Technology, Cambridge, MA 02139}
\email{clarkbar@math.mit.edu}

\author{D.M. Kan}
\address{Department of Mathematics, Massachusetts Institute of
  Technology, Cambridge, MA 02139}

\date{\today}

\maketitle

\begin{abstract}
  In this note we consider \emph{partial model categories}, by which
  we mean relative categories that satisfy a weakened version of the
  model category axioms involving only the weak equivalences.
  More precisely, a partial model category will be a \emph{relative
    category} that has the \emph{two out of six property} and admits a
  \emph{$3$-arrow calculus}.

  We then show that Charles Rezk's result that \emph{the simplicial
    space obtained from a simplicial model category by taking a
    Reedy fibrant replacement of its simplicial nerve is a
    complete Segal space} also holds for these partial model
  categories.

  We also note that conversely \emph{every complete Segal space is
    Reedy equivalent to the simplicial nerve of a partial model
    category} and in fact of a homotopically full subcategory of a
  category of diagrams of simplicial sets.
\end{abstract}

\section{Introduction}
\label{sec:Intro}

\subsection{Background and motivation}
\label{sec:Bckgrnd}

\begin{enumerate}
\item \label{BgRlCt} Given a \emph{relative category} $(\cat C,\cat
  U)$ (i.e.\ a pair consisting of a category $\cat C$ and a
  subcategory $\cat U$ which contains all the objects and of which the
  maps are called \emph{weak equivalences}), one can form the
  \emph{localization} of $\cat C$ with respect to $\cat U$, i.e.\ the
  category $\Ho(\cat C,\cat U)$ (called its \emph{homotopy category}),
  obtained from $\cat C$ by formally inverting the weak equivalences.

  However it was noted in \cite{DK1} that one can also form a
  \emph{simplicial localization} $L(\cat C,\cat U)$ of $\cat C$ with
  respect to $\cat U$, which is a simplicial category (i.e.\ a
  category enriched in simplicial sets) which has the same objects as
  $\cat C$ and which has the property that the category obtained by
  replacing each simplicial set by its set of components is exactly
  $\Ho(\cat C,\cat U)$.

  Moreover it was noted in \cites{DK2,DK3} that
  \begin{itemize}
  \item if $\cat M$ is a simplicial model category and $X$ and $Y$ are
    respectively cofibrant and fibrant objects of $\cat M$, then the
    function complex $\cat M_{*}(X,Y)$ has the same homotopy type as
    the simplicial set $L(\cat M,\cat W)(X,Y)$, where $\cat W \subset
    \cat M$ denotes its category of weak equivalences.
  \end{itemize}
  A key step in the proof of this result was the observation that if
  \begin{itemize}
  \item a relative category $(\cat C,\cat U)$ with the two out of
    three property admits a \emph{$3$-arrow calculus} (which means
    that there exists subcategories $\cat U_{c}$ and $\cat U_{f}
    \subset \cat U$ which have some of the properties of the
    categories of the trivial cofibrations and trivial fibrations in a
    model category)
  \end{itemize}
  then
  \begin{itemize}
  \item for every two objects $C_{1},C_{2}\in\cat C$, the homotopy
    type of $L(\cat C,\cat U)(C_{1},C_{2})$ admits a rather simple
    description in terms of $3$-arrow zigzags
    \begin{displaymath}
      C_{1} \longleftarrow \cdot \longrightarrow \cdot
      \longleftarrow C_{2}
    \end{displaymath}
    in which the outside maps are weak equivalences.
  \end{itemize}
\item While \cites{DK1,DK2,DK3} dealt with attempts to understand the
  function complexes in simplicial model categories, \cite{DHKS}
  turned its attention to the homotopy limit and colimit functors.  It
  was noted there that if
  \begin{itemize}
  \item a \emph{relative category} $(\cat C,\cat U)$, instead of the
    usual two out of three property, had the stronger \emph{two out of
      six property} that, for every three maps $r$, $s$ and $t \in
    \cat C$ for which $sr$ and $ts$ were defined, if the \emph{two}
    maps $sr$ and $ts$ were in $\cat U$, then so were the other
    \emph{four} maps $r$, $s$, $t$ and $tsr$,
  \end{itemize}
  then
  \begin{itemize}
  \item one could define homotopy limit and colimit functors which, if
    they existed, were homotopically unique, and give sufficient
    conditions for their existence and composability
  \end{itemize}
  and
  \begin{itemize}
  \item these sufficient conditions could be simplified if $(\cat
    C,\cat U)$ was \emph{saturated} (i.e.\ a map in $\cat C$ was a
    weak equivalence iff its image in $\Ho(\cat C,\cat U)$ was an
    isomorphism)
  \end{itemize}
  while
  \begin{itemize}
  \item a sufficient condition for such saturation was the presence of
    a \emph{$3$-arrow calculus} (i).
  \end{itemize}
\item All this suggests that the notion of a \emph{relative category}
  which has the \emph{two out of six property} and admits a
  \emph{$3$-arrow calculus} is a useful one and deserves further
  investigation.  As these three conditions are essentially parts of
  the model category axioms which only involve the weak equivalences
  we will refer to categories with weak equivalences with these three
  properties as \emph{partial model categories}.
\end{enumerate}

\subsection{Some examples of partial model categories}
\label{sec:ExWMC}

Some rather obvious examples of partial model categories are:
\begin{enumerate}
\item \label{ExMdCt} Model categories, and in particular the relative category $\cat S$ of simplicial sets.
\item \label{ExFlSbt} If $(\cat C,\cat U)$ is a partial model
  category, then so are all its \textbf{homotopically full} relative
  categories, i.e.\ all relative categories of the form $(\cat C',\cat
  C'\intersect\cat U)$ where $\cat C'$ is a full subcategory of $\cat
  C$ with the property that, for every object $C'\in\cat C'$, all
  objects of $\cat C$ which are weakly equivalent to $C'$ are also in
  $\cat C'$.
\item \label{ExRlFnc} If $(\cat C,\cat U)$ is a partial model
  category, then so is, for every relative category $(\cat A,\cat X)$
  the relative category $(\cat C,\cat U)^{\cat A,\cat X}$ of the
  relative functors $(\cat A,\cat X) \to (\cat C,\cat U)$.
\end{enumerate}

\intro
Our aim in this note now is to prove the following.
\subsection{A generalization of a result of Charles Rezk}
\label{sec:GenRezk}

\begin{enumerate}
\item \label{GenRezki} Rezk's result \cite{R}*{8.3} that ``for every
  \emph{simplicial model category} the Reedy fibrant replacement of
  its simplicial nerve is a \emph{complete Segal space}'' also holds
  for \emph{partial model categories}.
\end{enumerate}
Moreover, conversely,
\begin{resumeenumerate}{2}
\item \label{GenRezkii} every \emph{complete Segal space} is Reedy
  equivalent to the simplicial nerve of a \emph{partial model
    category}, and in fact of a homotopically full \eqref{ExFlSbt}
  relative subcategory of a category of diagrams of simplicial sets.
\end{resumeenumerate}

As if often the case with more general results, the proof of (i) is
simpler than Rezk's.  In proving the Segal part of his result, i.e.\
showing that certain fibre products (which are iterated pullbacks) are
homotopy fibre products (which are iterated homotopy pullbacks) he
relied heavily on the simplicial structure of his model category.
However it turns out that, in view of the \emph{Quillen Theorem
  B$_{3}$ for homotopy pullbacks} of \cite{BK3}*{8.2--4} this problem
reduces to a rather straightforward calculation which only involves
the partial model structure.

The proof of (ii) however requires a very different argument involving the following.

\subsection{A relative Yoneda embedding and partial modelization}
\label{sec:RlYnEmb}

Given a relative category $(\cat C,\cat U)$ we construct a
\emph{relative Yoneda embedding} (\ref{ExMdCt} and \ref{ExRlFnc})
\begin{displaymath}
  y\colon (\cat C,\cat U)\longrightarrow \cat S^{\cat C\op,\cat U\op}
\end{displaymath}
and note that
\begin{itemize}
\item the \textbf{essential image} $Ey$ of $y$, i.e.\ the
  homotopically full \eqref{ExFlSbt} relative subcategory of $\cat
  S^{\cat C\op,\cat U\op}$ (and hence of $\cat S^{\cat C\op}$) spanned
  by the image of $y$, is a partial model category
\end{itemize}
and that
\begin{itemize}
\item the inclusion $(\cat C,\cat U)\to Ey$ is a
  \textbf{DK-equivalence} i.e.\ its simplicial localization is a weak
  equivalence of simplicial categories \cite{Be} or equivalently
  \cite{BK2}*{1.8} its simplicial nerve is a Rezk (i.e.\ complete Segal)
  equivalence of simplicial spaces.
\end{itemize}

\subsection{Organization of the paper}
\label{sec:Org}

There are four more sections.

In \S\ref{sec:PrtMdCt} we introduce partial model categories and
discuss a few immediate consequences of their definition.  The notion
of a $3$-arrow calculus is slightly stronger than the one used in
\cites{DK2,DK3} and \cite{DHKS} in that in addition to the functorial
factorization we require that the trivial cofibration-like
subcategory be closed under pushouts and that the trivial fibration-like subcategory be closed under pullbacks.

In \S \ref{sec:GenRzk} and \S\ref{sec:RezkConverse} we then state and
prove the results that were mentioned in \ref{GenRezki} and
\ref{GenRezkii} respectively, while \S\ref{sec:PrtMdlLm} deals with the
partial modelization \eqref{sec:RlYnEmb}.

\section{Partial model categories}
\label{sec:PrtMdCt}

In this section we introduce \emph{partial model categories} and
discuss a few immediate consequences of their definition. 
\subsection{Partial model categories}
\label{sec:PrtModelCat}

A \textbf{partial model category} will be a pair $(\cat C,\cat W)$
consisting of a category $\cat C$ and a subcategory $\cat W\subset\cat
C$ (the maps of which will be called \textbf{weak equivalences})
which, roughly speaking, satisfies those parts of the model category
axioms (as for instance reformulated in \cite{DHKS}*{9.1}) which
involve only the weak equivalences.  More precisely we require that
\begin{lettered}
\item \label{reqRlCt} $(\cat C,\cat W)$ be a \textbf{relative
    category}, that $\cat W$ contains all the objects of $\cat C$ (and
  hence also their identity maps),
\item $(\cat C,\cat W)$ has the \textbf{two out of six property} that,
  if $r$, $s$ and $t$ are maps in $\cat C$ such that the \emph{two}
  compositions $sr$ and $ts$ exists and are in $\cat W$, then the
  \emph{four} maps $r$, $s$, $t$ and $tsr$ are also in $\cat W$ (which
  together with \ref{reqRlCt}) readily implies that \emph{$(\cat
    C,\cat W)$ has the two out of three property and that $\cat W$
    contains all the isomorphisms}.
\item $(\cat C,\cat W)$ admits a \textbf{$3$-arrow calculus}, i.e.\
  there exists subcategories $\cat U,\cat V\subset \cat W$ which
  behave very much like the categories of the trivial cofibrations and
  the trivial fibrations in a model category in the sense that
  \begin{enumerate}
  \item for every map $u \in \cat U$, its \emph{pushouts} in $\cat C$
    exist and are again in $\cat U$,
  \item for every map $v \in \cat V$, it' \emph{pullbacks} in $\cat C$
    exist and are again in $\cat V$, and
  \item \label{FuncFact} the maps $w \in \cat W$ admit a
    \emph{functorial factorization} $w = vu$ with $u \in \cat U$ and
    $v \in \cat V$ (which implies that $\cat U$ and $\cat V$ contain
    all the objects).
  \end{enumerate}
\end{lettered}

It should be noted that conditions (i) and (ii) are \emph{stronger}
than the ones that were used in \cite{DK2} and \cite{DHKS}.  However
we prefer them as they are cleaner and easier to work with and are
likely to be usually automatically satisfied.

\intro
One then readily verifies that the following are
\subsection{Examples of partial model categories}
\label{sec:ExPrtMC}

\begin{enumerate}
\item For every model category its underlying relative category is (of
  course) a partial model category.
\item Every homotopically full \eqref{ExFlSbt} relative subcategory of
  a partial model category is again a partial model category.
\item For every partial model category $(\cat C,\cat U)$ and relative
  category $(\cat A,\cat X)$, the relative functor category $(\cat
  C,\cat U)^{\cat A,\cat X}$ \eqref{ExRlFnc} is again a partial model category.
\item\label{item:relcatUU} For every partial model category $(\cat C,\cat U)$ the relative
  category $(\cat U,\cat U)$ is a partial model category.
\end{enumerate}

\section{A generalization of a result of Rezk}
\label{sec:GenRzk}

\subsection{Saturation \protect{\cite{DHKS}*{36.4}}}
\label{Prop:PrtMCSat}
\emph{Every partial model category $(\cat C,\cat W)$ is
  \textbf{saturated} in the sense that a map of $\cat C$ is in $\cat
  W$ iff it goes to an isomorphism in the \textbf{homotopy category}
  $\Ho(\cat C,\cat W)$, i.e.\ the category obtained from $\cat C$ by
  ``formally inverting'' the weak equivalences.}

\intro
Using this we will now show that a much stronger result of Rezk on
simplicial model categories \cite{R}*{8.3} also holds for partial model
categories.

Before formulating this we first recall
\subsection{Rezk's complete Segal model structure}
\label{sec:RzkCSgMC}

In \cite{R}
\begin{lettered}
\item Rezk constructed a ``homotopy theory of homotopy theories''
  model structure on the category $\simp\cat S$ of simplicial spaces
  (i.e.\ bisimplicial sets) by means of an appropriate left Bousfield
  localization of the Reedy model structure, the fibrant objects of
  which he referred to as \textbf{complete Segal spaces}, and
\item described a \textbf{Rezk} (or \textbf{simplicial}) \textbf{nerve
    functor} $N$ from the category $\RelCat$ of relative categories
  \eqref{sec:PrtModelCat} and relative functors between them to
  $\simp\cat S$ which sends a relative category $(\cat C,\cat W)$ to
  the simplicial space which in dimension $k \ge 0$ has as its
  $n$-simplices ($n \ge 0$) the commutative squares of the form
  \begin{displaymath}
    \xymatrix{
      {\cdot} \ar[r]^-{c_{1}} \ar[d]_{w_{1}}
      & {\cdot} \ar@{}[r]|\spacedcdots
      & {\cdot} \ar[r]^-{c_{k}}
      & {\cdot} \ar[d]\\
      {\cdot} \ar@{}[d]|{\vdots}
      &&& {\cdot} \ar@{}[d]|{\vdots}\\
      {\cdot} \ar[d]_{w_{n}}
      &&& {\cdot} \ar[d]\\
      {\cdot} \ar[r]
      & {\cdot} \ar@{}[r]|\spacedcdots
      & {\cdot} \ar[r]
      & {\cdot}
    }
  \end{displaymath}
  in which the vertical maps are in $\cat W$.
\end{lettered}
He then noted that
\begin{itemize}
\item[$*$] \emph{for every simplicial model category $\cat M$ one (and
    hence every) Reedy fibrant replacement of the simplicial space
    $N\cat M$ is a complete Segal space.}
\end{itemize}

Our aim thus is to prove the following.
\subsection{A generalization of Rezk's result}
\label{Thm:ReFibComSeg}

\emph{If $(\cat C,\cat W)$ is a partial model category
  \eqref{sec:PrtModelCat}, then one (and hence every) Reedy fibrant
  replacement of $N(\cat C,\cat W)$ is a complete Segal space.}

\emsection{Proof}
The proof consists of two parts, a \emph{Segal} part and a
\emph{completion} part.

To deal with the Segal part
\begin{enumerate}
\item let for every integer $k \ge 0$, $\cat A_{k}$ denote the
  category which has as its objects the sequences
  \begin{displaymath}
    \xymatrix{
        {\cdot} \ar[r]^{a_{1}}
        & {\cdot} \ar@{}[r]|{\spacedcdots}
        & {\cdot} \ar[r]^{a_{k}}
        & {\cdot}
    } \qquad\text{in $\cat C$}
  \end{displaymath}
  and as its maps the commutative diagrams of the form
  \begin{displaymath}
    \vcenter{\xymatrix{
        {\cdot} \ar[r]^{a_{1}} \ar[d]
        & {\cdot} \ar[d] \ar@{}[r]|{\spacedcdots}
        & {\cdot} \ar[r]^{a_{k}} \ar[d]
        & {\cdot} \ar[d]\\
        {\cdot} \ar[r]^{a'_{1}}
        & {\cdot} \ar@{}[r]|{\spacedcdots}
        & {\cdot} \ar[r]^{a'_{k}}
        & {\cdot}
      } 
    } \qquad\text{in $\cat C$}
  \end{displaymath}
  in which the vertical maps are in $\cat W$.
\end{enumerate}
Then we have to show that, for every integer $k \ge 2$, the pullback
square
\begin{displaymath}
  \xymatrix{
    {\cat A_{k}} \ar[r] \ar[d]
    & {\cat A_{k-1}} \ar[d]\\
    {\cat A_{1}} \ar[r]
    & *++[r]{\cat A_{0} = \cat W}
  }
\end{displaymath}
is a \emph{homotopy pullback square}.

To do this
\begin{resumeenumerate}{2}
\item for every integer $k \ge 2$, denote by $\cat B_{k}$ the
  category which has as its objects the zigzags
  \begin{displaymath}
    \xymatrix{
      {\cdot} \ar[r]^{b_{1}}
      & {\cdot} \ar[r]^{x}
      & {\cdot}
      & {\cdot} \ar[l]_{w} \ar[r]^{y}
      & {\cdot} \ar[r]^{b_{2}}
      & {\cdot} \ar@{}[r]|{\spacedcdots}
      & {\cdot} \ar[r]^{b_{k}}
      & {\cdot}
    } \qquad\text{in $\cat C$}
  \end{displaymath}
  in which $x$, $y$ and $w$ are in $\cat W$ and as its maps the
  commutative diagrams of the form
  \begin{displaymath}
    \vcenter{\xymatrix{
        {\cdot} \ar[r]^{b_{1}} \ar[d]
        & {\cdot} \ar[r]^{x} \ar[d]
        & {\cdot} \ar[d]
        & {\cdot} \ar[l]_{w} \ar[r]^{y} \ar[d]
        & {\cdot} \ar[r]^{b_{2}} \ar[d]
        & {\cdot} \ar@{}[r]|{\spacedcdots} \ar[d]
        & {\cdot} \ar[r]^{b_{k}} \ar[d]
        & {\cdot} \ar[d]\\
        {\cdot} \ar[r]^{b'_{1}}
        & {\cdot} \ar[r]^{x'}
        & {\cdot}
        & {\cdot} \ar[l]_{w'} \ar[r]^{y'}
        & {\cdot} \ar[r]^{b'_{2}}
        & {\cdot} \ar@{}[r]|{\spacedcdots}
        & {\cdot} \ar[r]^{b'_{k}}
        & {\cdot}
      } 
    } \qquad\text{in $\cat C$}
  \end{displaymath}
  in which the vertical maps are in $\cat W$, and
\item for every integer $k \ge 2$, denote by
  \begin{displaymath}
    h_{k}\colon \cat A_{k} \longrightarrow \cat B_{k}
  \end{displaymath}
  the monomorphism which between the first two maps inserts
  three identity maps, and denote by
  \begin{displaymath}
    \cat A'_{k}\subset \cat B_{k}
  \end{displaymath}
  the image of $\cat A_{k}$ under $h_{k}$.
\end{resumeenumerate}

In view of the Quillen Theorem B$_{3}$ for homotopy pullbacks
\cite{BK3}*{8.2--4} and the fact that, in view of \ref{item:relcatUU}, $\cat A_{0} = \cat W$ has property
$C_{3}$, it then suffices to show that, for every
integer $k \ge 2$, the inclusion $i\colon \cat A'_{k} \to \cat B_{k}$
is a \emph{homotopy equivalence}, i.e.\ that there exists a retraction
$r\colon \cat B_{k} \to \cat A'_{k}$ such that the compositions $ir$
and $ri$ are naturally weakly equivalent to the identity functor of
$\cat B_{k}$ and $\cat A'_{k}$ respectively.

Such a retraction, together with a zigzag of natural weak equivalences
connecting the functors $ir$ and $1_{\cat B_{k}}$ can be obtained by
means of the following (natural) commutative diagram in $\cat C$
\begin{displaymath}
  \xymatrix{
    {\cdot} \ar[r]^{b_{1}} \ar[d]
    & {\cdot} \ar[r]^{x} \ar[d]^{x}
    & {\cdot} \ar[d]
    & {\cdot} \ar[l]_{w} \ar[r]^{y} \ar[d]
    & {\cdot} \ar[r]^{b_{2}} \ar[d]
    & {\cdot} \ar@{}[r]|{\spacedcdots} \ar[d]
    & {\cdot} \ar[r]^{b_{k}} \ar[d]
    & {\cdot} \ar[d]\\
    {\cdot} \ar[r]^{xb_{1}}
    & {\cdot} \ar[r]
    & {\cdot}
    & {\cdot} \ar[l]_{w} \ar[r]^{y}
    & {\cdot} \ar[r]^{b_{2}}
    & {\cdot} \ar@{}[r]|{\spacedcdots}
    & {\cdot} \ar[r]^{b_{k}}
    & {\cdot}\\
    {\cdot} \ar[r]^{xb_{1}} \ar[u] \ar[d]
    & {\cdot} \ar[r]  \ar[u] \ar[d]
    & {\cdot}  \ar[u] \ar[d]
    & {\cdot} \ar[l]_{w} \ar[r]  \ar[u] \ar[d]_{u_{1}}
    & {\cdot} \ar[r]^{b_{2}y}  \ar[u]_{y} \ar[d]_{u_{1}}
    & {\cdot} \ar@{}[r]|{\spacedcdots}  \ar[u] \ar[d]_{u_{2}}
    & {\cdot} \ar[r]^{b_{k}}  \ar[u] \ar[d]^{u_{k-1}}
    & {\cdot}  \ar[u] \ar[d]^{u_{k}}\\
    {\cdot} \ar[r]^{xb_{1}}
    & {\cdot} \ar[r]
    & {\cdot}
    & {\cdot} \ar[l]_{v_{1}} \ar[r]
    & {\cdot} \ar[r]^{\overline{b_{2}y}}
    & {\cdot} \ar@{}[r]|{\spacedcdots}
    & {\cdot} \ar[r]^{\overline{b_{k}}}
    & {\cdot}\\
    {\cdot} \ar[r]^{\overline{xb_{1}}} \ar[u]^{\overline{v_{1}}}
    & {\cdot} \ar[r]  \ar[u]_{v_{1}}
    & {\cdot}  \ar[u]_{v_{1}}
    & {\cdot} \ar[l] \ar[r]  \ar[u]
    & {\cdot} \ar[r]^{\overline{b_{2}y}}  \ar[u]
    & {\cdot} \ar@{}[r]|{\spacedcdots}  \ar[u]
    & {\cdot} \ar[r]^{\overline{b_{k}}}  \ar[u]
    & {\cdot} \ar[u]
  }
\end{displaymath}
in which all the unmarked arrows are identity maps, $w = v_{1}u_{1}$
with $u_{1}\in\cat U$ and $v_{1}\in\cat V$ (\ref{sec:PrtModelCat}c)
and the squares involving two $u$'s are pushout squares and those
involving two $v$'s are pullback squares.

On $\cat A'_{k}$ this zigzag reduces to the zigzag
\begin{displaymath}
  \xymatrix{
    {\cdot} \ar[r]^{b_{1}} \ar[d]
    & {\cdot} \ar[r] \ar[d]
    & {\cdot} \ar[d]
    & {\cdot} \ar[l] \ar[r] \ar[d]_{u_{1}}
    & {\cdot} \ar[r]^{b_{2}} \ar[d]_{u_{1}}
    & {\cdot} \ar@{}[r]|{\spacedcdots} \ar[d]_{u_{2}}
    & {\cdot} \ar[r]^{b_{k}} \ar[d]_{u_{k-1}}
    & {\cdot} \ar[d]_{u_{k}}\\
    {\cdot} \ar[r]^{b_{1}}
    & {\cdot} \ar[r]
    & {\cdot}
    & {\cdot} \ar[l]_{v_{1}} \ar[r]
    & {\cdot} \ar[r]^{\bar b_{2}}
    & {\cdot} \ar@{}[r]|{\spacedcdots}
    & {\cdot} \ar[r]^{\bar b_{k}}
    & {\cdot}\\
    {\cdot} \ar[r]^{\bar b_{1}} \ar[u]^{\bar v_{1}}
    & {\cdot} \ar[r]  \ar[u]_{v_{1}}
    & {\cdot}  \ar[u]_{v_{1}}
    & {\cdot} \ar[l] \ar[r]  \ar[u]
    & {\cdot} \ar[r]^{\bar b_{2}}  \ar[u]
    & {\cdot} \ar@{}[r]|{\spacedcdots}  \ar[u]
    & {\cdot} \ar[r]^{\bar b_{k}}  \ar[u]
    & {\cdot}  \ar[u]
  }
\end{displaymath}
which does not completely lie inside $\cat A'_{k}$.  To remedy this,
i.e.\ to get a natural weak equivalence connecting the top with the
bottom inside $\cat A'_{k}$ we note the existence of the zigzag
\begin{displaymath}
  \xymatrix{
    {\cdot} \ar[r]^{b_{1}} \ar[d]
    & {\cdot} \ar[r] \ar[d]
    & {\cdot} \ar[d]
    & {\cdot} \ar[l] \ar[r] \ar[d]_{u_{1}}
    & {\cdot} \ar[r]^{b_{2}} \ar[d]_{u_{1}}
    & {\cdot} \ar@{}[r]|{\spacedcdots} \ar[d]_{u_{2}}
    & {\cdot} \ar[r]^{b_{k}} \ar[d]_{u_{k-1}}
    & {\cdot} \ar[d]_{u_{k}}\\
    {\cdot} \ar[r]^{b_{1}} \ar[d]
    & {\cdot} \ar[r] \ar[d]
    & {\cdot}  \ar[d]
    & {\cdot} \ar[l]_{v_{1}} \ar[r] \ar[d]_{v_{1}}
    & {\cdot} \ar[r]^{\bar b_{2}} \ar[d]_{v_{1}}
    & {\cdot} \ar@{}[r]|{\spacedcdots} \ar[d]_{v_{2}}
    & {\cdot} \ar[r]^{\bar b_{k}} \ar[d]_{v_{k-1}}
    & {\cdot} \ar[d]_{v_{k}}\\
    {\cdot} \ar[r]^{b_{1}}
    & {\cdot} \ar[r]
    & {\cdot} 
    & {\cdot} \ar[l] \ar[r]
    & {\cdot} \ar[r]^{b_{2}}
    & {\cdot} \ar@{}[r]|{\spacedcdots}
    & {\cdot} \ar[r]^{b_{k}}
    & {\cdot} 
  }
\end{displaymath}
in which the bottom row is obtained from the top row by pushing out
along $v_{1}u_{1}$ which is an identity map.  Combining the bottom
halves of the last two diagrams we now get two composable natural weak
equivalences
\begin{displaymath}
  \xymatrix{
    {\cdot} \ar[r]^{b_{1}}
    & {\cdot} \ar[r]
    & {\cdot}
    & {\cdot} \ar[l] \ar[r]
    & {\cdot} \ar[r]^{b_{2}}
    & {\cdot} \ar@{}[r]|{\spacedcdots}
    & {\cdot} \ar[r]^{b_{k}}
    & {\cdot}\\
    {\cdot} \ar[r]^{b_{1}} \ar[u]
    & {\cdot} \ar[r] \ar[u]
    & {\cdot} \ar[u]
    & {\cdot} \ar[l]_{v_{1}} \ar[r] \ar[u]^{v_{1}}
    & {\cdot} \ar[r]^{\bar b_{2}} \ar[u]^{v_{1}}
    & {\cdot} \ar@{}[r]|{\spacedcdots} \ar[u]^{v_{2}}
    & {\cdot} \ar[r]^{\bar b_{k}} \ar[u]^{v_{k-1}}
    & {\cdot} \ar[u]^{v_{k}}\\
    {\cdot} \ar[r]^{\bar b_{1}} \ar[u]^{\bar v_{1}}
    & {\cdot} \ar[r] \ar[u]_{v_{1}}
    & {\cdot} \ar[u]_{v_{1}}
    & {\cdot} \ar[l] \ar[r] \ar[u]
    & {\cdot} \ar[r]^{\bar b_{2}} \ar[u]
    & {\cdot} \ar@{}[r]|{\spacedcdots} \ar[u]
    & {\cdot} \ar[r]^{\bar b_{k}} \ar[u]
    & {\cdot} \ar[u]
  }
\end{displaymath}
of which the composition
\begin{displaymath}
  \xymatrix{
    {\cdot} \ar[r]^{b_{1}}
    & {\cdot} \ar[r]
    & {\cdot}
    & {\cdot} \ar[l] \ar[r]
    & {\cdot} \ar[r]^{b_{2}}
    & {\cdot} \ar@{}[r]|{\spacedcdots}
    & {\cdot} \ar[r]^{b_{k}}
    & {\cdot}\\
    {\cdot} \ar[r]^{\bar b_{1}} \ar[u]^{\bar v_{1}}
    & {\cdot} \ar[r] \ar[u]_{v}
    & {\cdot} \ar[u]_{u}
    & {\cdot} \ar[l] \ar[r] \ar[u]^{v_{1}}
    & {\cdot} \ar[r]^{\bar b_{2}} \ar[u]^{v_{1}}
    & {\cdot} \ar@{}[r]|{\spacedcdots} \ar[u]^{v_{2}}
    & {\cdot} \ar[r]^{\bar b_{k}} \ar[u]^{v_{k-1}}
    & {\cdot} \ar[u]^{v_{k}}
  }
\end{displaymath}
yields the desired natural weak equivalence between $ri$ and $1_{\cat
  A'}$.

It thus remains to deal with the completeness part of the proof.
However this is essentially the same as Rezk's proof of \cite{R}*{8.3}
in view of the fact that the partial model category $(\cat C,\cat W)$
is saturated \eqref{Prop:PrtMCSat}.

\section{A converse of Rezk's result}
\label{sec:RezkConverse}

We now prove

\subsection{A converse of Rezk's result}
\label{thm:CmpSegEq}
\begin{em}
  Every complete Segal space is Reedy equivalent to the simplicial
  nerve of a partial model category and in fact of a homotopically
  full relative subcategory of a category of diagrams of simplicial
  sets.
\end{em}

\smallskip

The key to this is a \emph{partial modelization lemma} which we will
state in \ref{sec:PartMdlLm} but prove in \S \ref{sec:PrtMdlLm} below.
Its formulation requires the following.
\subsection{A relative Yoneda embedding}
\label{sec:RelYonEmb}

Let $\LH$ denote the hammock localization of \cite{DK2}.  Given a
relative category $(\cat C,\cat U)$, its \textbf{relative Yoneda
  embedding} will be the relative functor between relative categories
\begin{displaymath}
  y = y_{\cat C,\cat U}\colon (\cat C,\cat U) \longrightarrow
  \cat S^{\cat C\op,\cat U\op}
\end{displaymath}
which sends each object $A \in \cat C$ to the relative functor
$yA\colon (\cat C\op,\cat U\op) \to \cat S$ which sends each object
$B \in \cat C\op$ to the simplicial set $\LH(\cat C,\cat U)(B,A)$.

\subsection{The partial modelization}
\label{sec:PartMdlLm}

\begin{em}
  Given a relative category $(\cat C,\cat U)$,
  \begin{enumerate}
  \item \label{PartMdlLmi} the essential image $Ey$
    \eqref{sec:RlYnEmb} of its Yoneda embedding\eqref{sec:RelYonEmb}
    is a partial model category and in fact a homotopically full
    \eqref{ExFlSbt} relative subcategory of a category of diagrams of
    simplicial sets, and
  \item \label{PartMdlLmii} the embedding $e\colon (\cat C,\cat U) \to
    Ey$ is a DK-equivalence \eqref{sec:RlYnEmb}.
  \end{enumerate}
\end{em}

\intro
Using this we now can give
\subsection{A proof of \ref{thm:CmpSegEq}}
\label{sec:Prf4p1}

First recall from \cite{BK1}*{5.3 and 4.4} the existence of
\begin{enumerate}
\item \label{prf4p1i} an adjunction $\adj{K_{\xi}}{\RelCat}{\simp\cat
    S}{N_{\xi}}$ of which the unit $\eta\colon 1 \to N_{\xi}K_{\xi}$
  is a natural Reedy equivalence, and
\item \label{prf4p1ii} a natural Reedy equivalence $\pi^{*}\colon N
  \to N_{\xi}$ (\ref{sec:RzkCSgMC}b)
\end{enumerate}
and from \cite{R}*{7.2} that
\begin{resumeenumerate}{3}
\item \label{prf4p1iii} every Reedy equivalence in $\simp\cat S$ is a
  Rezk equivalence \eqref{PartMdlLmii} and every Rezk equivalence
  between two complete Segal spaces is a Reedy equivalence.
\end{resumeenumerate}

Given a complete Segal space $X$ one then can consider the zigzag
\eqref{sec:GenRezk}
\begin{displaymath}
  \xymatrix{
    {X} \ar[r]^-{\eta}
    & {N_{\xi}K_{\xi}X}
    & {NK_{\xi}X} \ar[r]^-{e} \ar[l]_-{\pi^{*}}
    & {NEy_{K_{\xi},X}}
  }
\end{displaymath}
in which, in view of (i) and (ii) above and \ref{PartMdlLmii}
respectively, the first two maps are Reedy equivalences, while the
third is a Rezk equivalence, and note that it follows from (iii) above
and \ref{Thm:ReFibComSeg} that every Reedy fibrant replacement of the
partial model category $Ey^{}_{K_{\xi}X}$ \eqref{PartMdlLmi} is Reedy
equivalent to $X$.

\section{A proof of the partial modelization lemma \eqref{sec:PartMdlLm}}
\label{sec:PrtMdlLm}

In preparation for the proof of lemma~\ref{sec:PartMdlLm} (in
\ref{sec:PrfRlMdl} below) we first
\begin{itemize}
\item discuss in \ref{sec:RelPrtSmpCt} relative \emph{simplicial}
  categories and in particular relative \emph{partly simplicial} ones
  in which the weak equivalences form an ordinary category, and
\item review in \ref{sec:FullFaith} and \ref{sec:EssIm} the notions of
  \emph{fully faithfulness} and \emph{essential surjectivity} and of
  \emph{essential image} in the categories of categories, simplicial
  categories, relative categories and relative simplicial categories.
\end{itemize}

\subsection{Relative (partly) simplicial categories}
\label{sec:RelPrtSmpCt}

Let $\SCat$ denote the category of \textbf{simplicial categories},
i.e.\ categories enriched over simplicial sets, and let $\RelSCat$
denote the resulting category of \textbf{relative simplicial
  categories}, i.e.\ pairs consisting of a simplicial category and a
sub-simplicial category (of which the maps are called \textbf{weak
  equivalences}) that contains all the objects.  Then it turns out
that, for our purposes here, it is convenient to work in the somewhat
simpler full subcategory
\begin{displaymath}
  \RelPSCat\subset \RelSCat
\end{displaymath}
spanned by what we will call the \textbf{relative partly simplicial
  categories}, i.e.\ the objects of which \emph{the weak equivalences
  form an ordinary category}.

A \emph{simplicial} model category then can be considered as
\begin{itemize}
\item an object of $\RelCat$ consisting of the underlying \emph{model
    category} and its weak equivalences
\end{itemize}
or as
\begin{itemize}
\item an object of $\RelPSCat$ consisting of the larger
  \emph{simplicially enriched model category} and those same weak
  equivalences.
\end{itemize}

Moreover in the remainder of this paper \textbf{we will consider the
  category $\cat S$ of simplicial sets only as an object of
  $\RelPSCat$}.

An object $\cat L \in \SCat$ thus gives rise to
\begin{itemize}
\item an object $(\cat S^{\cat L},\sim) \in \RelCat$ in which $\cat
  S^{\cat L}$ denotes the (model) category which has as objects the
  simplicial functors $\cat L \to \cat S$ and as maps the natural
  transformations between them and $\sim$ denotes the subcategory of
  the natural weak equivalences, and
\item an object $(\cat S_{\ast}^{\cat L},\sim)\in\RelSCat$ in which $\cat S_{\ast}^{\cat L}$
  denotes the simplicial (model) category of the simplicial functors
  $\cat L \to \cat S$ (\cite{DK4}*{1.3(v)} and \cite{GJ}*{IX, 1.4})
  and $\sim$ is as above.
\end{itemize}

We end with noting that similarly an object $(\cat L,\cat Z) \in
\RelPSCat$ gives rise to
\begin{itemize}
\item an object $(\cat S^{\cat L,\cat Z},\sim) \in\RelCat$ which is
  the subobject of $(\cat S^{\cat L},\sim)$ spanned by the relative
  simplicial functors $(\cat L,\cat Z) \to \cat S$
\end{itemize}
and that
\begin{em}
  \begin{enumerate}
  \item \label{RelPrtSmpCti} if $\cat Z$ is \textbf{neglectible} in
    $\cat L$ in the sense that every map in $\cat Z$ goes to an
    isomorphism in $\Ho\cat L$, then $(\cat S^{\cat L,\cat Z},\sim) =
    (\cat S^{\cat L},\sim)$, and
  \item \label{RelPrtSmpCtii} for every object $(\cat C,\cat U) \in
    \RelCat$, the object $(\cat S^{\cat C\op,\cat U\op},\sim) \in
    \RelCat$ is exactly the same as the object $\cat S^{\cat C\op,\cat
      U\op}$ mentioned in \ref{sec:RelYonEmb}.
  \end{enumerate}
\end{em}

\subsection{Fully faithfulness and essential surjectivity}
\label{sec:FullFaith}

We will denote by $\LH$ not only the functor $\RelCat \to \SCat$ which
sends each object to its hammock localization \cite{DK2}*{2.1}, but
also the functor $\RelSCat \to \SCat$ which sends each object to the
diagonal of the bisimplicial category obtained from it by
dimensionwise application of the hammock localization
\cite{DK2}*{2.5}.

Then we recall the following.

A functor $f\colon \cat G \to \cat H$ between categories (respectively, simplicial categories)
is called \textbf{fully faithful} if, for every two objects
$G_{1},G_{2} \in \cat G$, it induces an isomorphism (resp.\ weak
equivalence) $\cat G(G_{1},G_{2}) \to \cat H(fG_{1},fG_{2})$, and
similarly a relative functor $f\colon (\cat C,\cat U) \to (\cat D,\cat
V)$ between relative categories (resp.\ relative simplicial categories) is called \textbf{fully
  faithful} if, for every two objects $C_{1},C_{2} \in \cat C$, it
induces a weak equivalence
\begin{displaymath}
  \LH(\cat C,\cat U)(C_{1},C_{2}) \longrightarrow
  \LH(\cat D,\cat V)(fC_{1},fC_{2}) \in \cat S
\end{displaymath}
which implies that
\begin{em}
  \begin{enumerate}
  \item \label{FullFaithi} if $f$ and $g$ are (relative) functors such
    that $gf$ is defined and $g$ is fully faithful, then $gf$ is fully
    faithful iff $f$ is so.
  \end{enumerate}
\end{em}

A functor $f\colon \cat G \to \cat H$ between categories (respectively, simplicial categories)
is called \textbf{essentially surjective} if every object in $\cat H$
is isomorphic in $\cat H$ (resp.\ $\Ho\cat H$) to an object in the
image of $f$ (resp.\ $\Ho f$), and similarly a relative functor
$f\colon (\cat C,\cat U) \to (\cat D,\cat V)$ between relative categories (resp.\ relative simplicial categories) is called \textbf{essentially surjective} if
the induced functor
\begin{displaymath}
  \LH f\colon \LH(\cat C,\cat U) \longrightarrow
  \LH(\cat D,\cat V)
\end{displaymath}
is so, which implies that
\begin{em}
  \begin{resumeenumerate}{2}
  \item \label{FullFaithii} if $f$ and $g$ are (relative) functors
    such that $gf$ is defined and $f$ is essentially surjective, then
    $gf$ is essentially surjective iff $g$ is so.
  \end{resumeenumerate}
\end{em}
Then
\begin{em}
  \begin{resumeenumerate}{3}
  \item \label{FullFaithiii} a map in $\Cat$ is an equivalence of
    categories iff it is fully faithful and essentially surjective,
    and
  \item \label{FullFaithiv} a map in $\RelCat$, $\SCat$ or $\RelSCat$
    is a DK-equivalence iff it is fully faithful and essentially
    surjective.
  \end{resumeenumerate}
\end{em}

\subsection{Essential images}
\label{sec:EssIm}

The \textbf{essential image} $Ef$ of a functor $f\colon \cat G \to
\cat H$ between categories (respectively, simplicial categories) is the full
subcategory (resp.\ full simplicial subcategory) of $\cat H$ spanned by the objects which are isomorphic in
$\cat H$ (resp.\ $\Ho\cat H$) to objects in the image of $f$ (resp.\ $\Ho f$) and
similarly the \textbf{essential image} $Ef$ of a relative functor
$f\colon (\cat C,\cat U) \to (\cat D,\cat V)$ between relative
categories (resp.\ relative simplicial categories) is defined by the pullback diagram
\begin{displaymath}
  \xymatrix{
    {Ef} \ar[r] \ar[d]
    & {(\cat D,\cat V)} \ar[d]\\
    {E\LH f} \ar[r]
    & {\LH(\cat D, \cat V)}
  }
\end{displaymath}
which implies that
\begin{em}
  \begin{enumerate}
  \item \label{EssImi} the resulting maps
    \begin{displaymath}
      \cat G \longrightarrow Ef
      \qquad\text{and}\qquad
      (\cat C,\cat U) \longrightarrow Ef
    \end{displaymath}
    and
    \begin{displaymath}
      Ef \longrightarrow \cat H
      \qquad\text{and}\qquad
      Ef \longrightarrow (\cat D,\cat V)
    \end{displaymath}
    are respectively essentially surjective and fully faithful.
  \end{enumerate}
\end{em}

We end with noting that
\begin{em}
  \begin{resumeenumerate}{2}
  \item \label{EssImii} the essential image defined in
    \ref{sec:RlYnEmb} is a special case of the ones defined above.
  \end{resumeenumerate}
\end{em}

\intro
Now we are finally ready for

\subsection{A proof of the relative modelization lemma
  \eqref{sec:PartMdlLm}}
\label{sec:PrfRlMdl}

It follows from \ref{EssImi} and (ii) that the map $e\colon (\cat
C,\cat U) \to Ef$ is essentially surjective and it thus
\eqref{FullFaithiv} remains to prove that it is also fully faithful.
To do this it suffices, in view of \ref{EssImii} and \ref{FullFaithi},
to show that, in the notation of \ref{RelPrtSmpCtii},
\begin{enumerate}
\item \label{PrfRlMdli} the Yoneda embedding $y\colon (\cat C,\cat U)
  \to (\cat S^{\cat C\op,\cat U\op},\sim)$ is fully faithful.
\end{enumerate}

For this we note that $y$ admits a factorization
\eqref{sec:RelPrtSmpCt}
\begin{displaymath}
  (\cat C,\cat U) \xrightarrow{y'}
  (\cat S^{\LH(\cat C\op,\cat U\op)},\sim)
  \overset{\ref{RelPrtSmpCti}}{=}
  (\cat S^{\LH(\cat C\op,\cat U\op),\cat U\op},\sim)
  \xrightarrow{(c\op)^{*}}
  (\cat S^{\cat C\op,\cat U\op},\sim)
\end{displaymath}
in which $y'$ sends each object $A \in \cat C$ to the simplicial
functor $\LH(\cat C\op,\cat U\op) \to \cat S$ which sends each object
$B \in \LH(\cat C\op,\cat U\op)$ to $\LH(B,A) \in \cat S$, and
$c\colon (\cat C\op,\cat U\op) \to \LH(\cat C\op,\cat U\op),\cat U\op$
is the obvious inclusion \cite{DK2}*{3.1}.  The latter is a
DK-equivalence \cite{BK2}*{3.2} and hence \cite{DK4}*{2.2} so is the
map $(c\op)^{*}$.  Hence, in view of \ref{FullFaithiv} and
\ref{FullFaithi}, the condition (i) above is equivalent to condition
\begin{resumeenumerate}{2}
\item \label{PrfRlMdlii} the map $y'\colon (\cat C,\cat U) \to (\cat
  S^{\LH(\cat C\op,\cat U\op)},\sim)$ is fully faithful.
\end{resumeenumerate}

To prove this we embed this map in the commutative diagram
\begin{displaymath}
  \xymatrix{
    {(\cat C,\cat U)} \ar[r]^-{y'} \ar[d]_{c}
    & {(\cat S^{\LH(\cat C\op,\cat U\op)},\sim)}
    \ar[d]^{\text{incl.}}\\
    {\bigl(\LH(\cat C,\cat U),\cat U\bigr)} \ar[r]^-{r'}
    & (\cat S_{*}^{\LH(\cat C\op,\cat U\op)},\sim)
  }
\end{displaymath}
in which the map in the right is as in \ref{sec:RelPrtSmpCt} and $r'$
is induced by the \emph{simplicial Yoneda embedding} of
\cite{DK4}*{1.3(vi)}
\begin{displaymath}
  r\colon \LH(\cat C,\cat U) \longrightarrow
  \cat S_{*}^{\LH(\cat C\op,\cat U\op)} \in \SCat
\end{displaymath}
which sends each object $A \in \LH(\cat C,\cat U)$ to the simplicial
functor $\LH(\cat C\op,\cat U\op) \to \cat S$ which sends each object
$B \in \LH(\cat C\op,\cat U\op)$ to $\LH(B,A) \in \cat S$.  The map
on the left is (see above) a DK-equivalence and so is, in view of
\cite{DK3}*{4.8} the map on the right and hence, to prove (ii), it
suffices to show that the bottom map is fully faithful.

For this we embed this map in the following diagram
\begin{displaymath}
  \xymatrix@C=4em{
    {\LH\bigl(\LH(\cat C,\cat U),\cat U\bigr)} \ar[r]^-{\LH r'}
    & {\LH(\cat S_{*}^{\LH(\cat C\op,\cat U\op)},\sim)}\\
    {\LH(\cat C,\cat U)} \ar[u] \ar[r]^-{r}
    & {\cat S_{*}^{\LH(\cat C\op,\cat U\op)}} \ar[u]
  }
\end{displaymath}
in which the vertical maps are the obvious inclusions
\cite{DK2}*{3.1}.  As both categories of weak equivalences are
neglectible \ref{RelPrtSmpCti}, it follows from \cite{DK1}*{6.4} that
both vertical maps are DK-equivalences.  Moreover it was noted in
\cite{DK4}*{1.3(vi)} that the map $r$ is fully faithful (and in fact
so in the strong sense that the required weak equivalences are
actually isomorphisms).  All this implies that $\LH r'$ is fully
faithful and so is therefore the map $r'$ itself.


\end{document}